\documentclass[14pt,twoside]{article}
\usepackage{graphicx}
\usepackage{amsmath,amsbsy,amsfonts}
\usepackage{color}
\def\dis{\displaystyle}
\def\nd{\noindent}

\newtheorem{lemma}{Lemma}[section]
\newtheorem{prop}{Proposition}[section]
\newtheorem{theorem}{Theorem}[section]

\newtheorem{coro}{Corollary}[section]

\newtheorem{definition}[theorem]{Definition}
\newcommand{\fim}{\hfill\rule{2mm}{2mm}}

\newcommand{\ds}{\displaystyle}
\def\dis{\displaystyle}

\let\Section=\section
\def\section{\setcounter{equation}{0}\Section}

\begin{document}
\title{ On a nonlocal multivalued problem  in an Orlicz-Sobolev space via Krasnoselskii's genus }

\author{ Giovany M. Figueiredo\thanks{Partially supported by CNPq/PQ
301242/2011-9 and 200237/2012-8 }  \\
\noindent Universidade Federal do Par\'a, Faculdade de Matem\'atica, \\
\noindent CEP: 66075-110, Bel\'em - Pa, Brazil \\
\noindent e-mail: giovany@ufpa.br  \vspace{0.5cm}\\
\noindent Jefferson A. Santos \thanks{Partially supported by CNPq-Brazil grant Casadinho/Procad 552.464/2011-2 } \\
\noindent Universidade Federal de Campina Grande,\\
\noindent Unidade Acad\^emica de Matem\'atica e Estat\'istica,\\
\noindent CEP:58109-970, Campina Grande - PB, Brazil\\
\noindent e-mail: jefferson@dme.ufcg.edu.br\vspace{0.5cm}\\
}\vspace{0.5cm}
\date{}

\pretolerance10000

\maketitle

\begin{abstract}
{This paper is concerned with the multiplicity of nontrivial
solutions in an Orlicz-Sobolev space for a nonlocal problem
involving N-functions and theory of locally Lispchitz continuous
functionals. More precisely, in this paper, we study a result of
multiplicity to the following multivalued elliptic problem:
$$
\left \{ \begin{array}{l}
-M\left(\displaystyle\int_\Omega \Phi(\mid\nabla u\mid)dx\right)div\big(\phi(\mid\nabla u\mid)\nabla u\big) -\phi(|u|)u\in  \partial F(u)  \ \mbox{in}\ \Omega,\\
u\in W_0^1L_\Phi(\Omega),                                    
\end{array}
\right.
$$
where  $\Omega\subset\mathbb{R}^{N}$ is a bounded smooth domain,
$N\geq 2$, $M$ is continuous function, $\Phi$ is an N-function with
$\Phi(t)=\displaystyle\int^{|t|}_{0}\phi(s)s \ ds$ and $\partial
F(t)$ is a generalized gradient of $F(t)$. We use genus theory to
obtain the main result.}
\end{abstract}

\maketitle

\section{Introduction}

The purpose of this article is investigate the multiplicity of
nontrivial solutions to the multivalued elliptic problem
$$
\left \{ \begin{array}{l}
-M\left(\displaystyle\int_\Omega \Phi(\mid\nabla u\mid)dx\right)div\big(\phi(\mid\nabla u\mid)\nabla u\big)-\phi(|u|)u\in  \partial F(u)  \ \mbox{in}\ \Omega,\\
u\in W_0^1L_\Phi(\Omega),                                    
\end{array}
\right.\leqno{(P)}
$$
 \noindent where $\Omega\subset\mathbb{R}^{N}$ is a bounded
smooth domain with $N\geq 2$, $F(t)= \displaystyle\int_0^t f(s)ds$
and
$$
\partial F(t)=\left\{s\in\mathbb{R}; F^0(t;r)\geq sr, \ r\in \mathbb{R}\right\}.
$$
Here $F^0(t;r)$ denotes the generalized directional derivative of
$t\mapsto F(t)$ in direction of $r$, that is,
$$
F^0(t;r)=\displaystyle \limsup_{h\rightarrow
t,s\downarrow0}\frac{F(h+sr)-F(h)}{s}.
$$
We shall assume in this work that $f(t)$ is locally bounded in
$\mathbb{R}$ and
$$
\underline{f}(t)=\displaystyle \lim_{\epsilon\downarrow 0}\mbox{ess
inf}\left\{ f(s);|s-t|<\epsilon\right\} \mbox{ and }
\overline{f}(t)=\lim_{\epsilon\downarrow 0}\mbox{ess sup}\left\{
f(s);|s-t|<\epsilon\right\} .
$$
It is well known that
$$
\partial F(t)=[\underline{f}(t),\overline{f}(t)], \mbox{ (see
\cite{chang})},
$$
and that, if $f(t)$ is continuous then $\partial F(t)=\{f(t)\}$.

Problem $(P)$  with $\phi(t)=2$, that is,
$$
\left \{ \begin{array}{l}
-M\left(\displaystyle\int_\Omega \mid\nabla u\mid^{2} dx\right)\Delta u - u \in \partial F(u) \ \mbox{in}\ \Omega,\\
u\in H^{1}_{0}(\Omega)                                    
\end{array}
\right.\leqno{(*)}
$$
is called nonlocal because of the presence of the term
$M\left(\dis\int_{\Omega}|\nabla u|^{2} dx \right)$ which implies
that the equation $(*)$ is no longer a pointwise identity.

The reader may consult $\cite{alvescorrea}$, $\cite{alvescorreama}$,
$\cite{GJ}$ and the references therein, for more
information on nonlocal problems.

On the other hand, in this study, the nonlinearity $f$ can be
discontinuous. There is by now an extensive literature on multivalued equations and we refer the reader to \cite{Goncalves},  \cite{Nascimento},  \cite{AlvesNascimento},  \cite{Santos},  \cite{Carvalho},  and references therein. The interest in the study of nonlinear partial
differential equations with discontinuous nonlinearities has
increased because many free boundary problems arising in
mathematical physics may be stated in this form.

Among these problems, we have the obstacle problem, the seepage surface problem,
and the Elenbaas equation, see for example \cite{chang},
\cite{chang1} and \cite{chang2}.

For enunciate the main result, we need to give some hypotheses on
the functions  $M, \phi$ and $f$.

The hypotheses on the function $\phi:\mathbb{R}^{+}\rightarrow
\mathbb{R}^{+}$ of $C^{1}$ class are the following:

\begin{description}

\item[($\phi_1$)] For all $t>0$,
$$
\phi(t)>0 \ \ \mbox{and} \ \ (\phi(t)t)'>0.
$$

\item[($\phi_2$)] There exist $l, m \in (1,N)$,  $l\leq m < l^{*}=
\displaystyle\frac{lN}{N-l}$ such that
$$
l\leq \frac{\phi(t)t^{2}}{\Phi(t)}\leq m,
$$
for $t> 0$, where $\Phi(t)=\displaystyle\int^{|t|}_{0}\phi(s)s ds$.
\end{description}

The hypothesis on the continuous function
$M:\mathbb{R}^{+}\rightarrow \mathbb{R}^{+}$ is the following:

\begin{description}

\item[($M_1$)]
There exist $k_{0},k_{1},\alpha, q_{0},q_{1}>0$ and
$b:\mathbb{R}\rightarrow \mathbb{R}$ of $C^{1}$ class such that
$$
k_0t^{\alpha}\leq M(t)\leq k_1t^{\alpha},
$$
$\alpha >\frac{q_1}{l}$, where
$$
m< q_0\leq \frac{b(t)t^2}{B(t)}\leq q_1 <l^*,
$$
for all $t>0$ with
$$
(b(t)t)'>0, \ t>0
$$
and
$$
B(t)=\int_0^tb(s)sds.
$$
\end{description}

The hypotheses on the function $f:\mathbb{R}\rightarrow \mathbb{R}$
are the following:
\begin{description}

\item[($f_1$)] For all $t \in \mathbb{R}$,

$$
f(t)=-f(-t).
$$

\item[($f_2$)] There exist $b_{0},b_{1}>0$ and $a_0\geq 0$ such that
$$
b_0b(t)t\leq f(t)\leq b_1b(t)t,\ |t|\geq a_0.
$$

\item[($f_3$)] There exists $a_0\geq 0$ such that
$$
f(t)=0,\ |t|\leq a_0.
$$
\end{description}

The main result of this paper is:

\begin{theorem}\label{teorema1}
Assume that conditions $(\phi_{1})$, $(\phi_{2})$, $(M_{1})$,
$(f_{1})-(f_{3})$ hold. Then for $a_0>0$ sufficiently small (or $a_0=0$), the problem $(P)$  has infinitely many solutions.
\end{theorem}

Below we show two graphs of functions that satisfy the hypotheses
$(f_{1})-(f_{3})$. Note that the second graph corresponds to a
function that has an enumerable number of points of discontinuity.
\begin{center}
\includegraphics[width=0.4\textwidth]{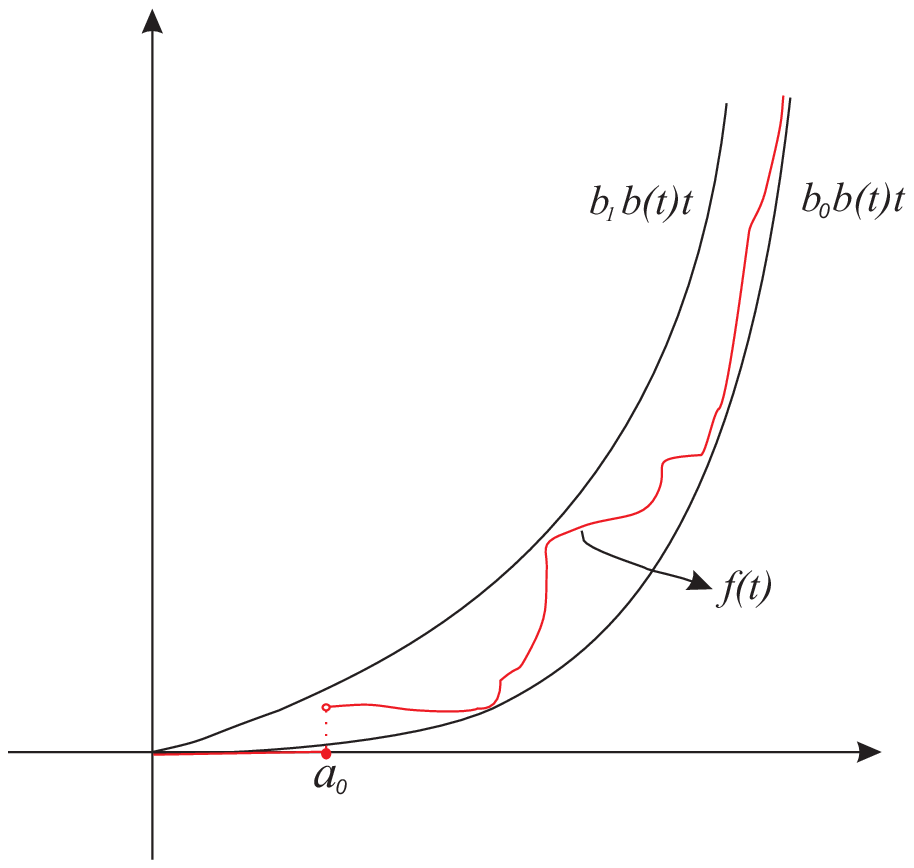}
\includegraphics[width=0.4\textwidth]{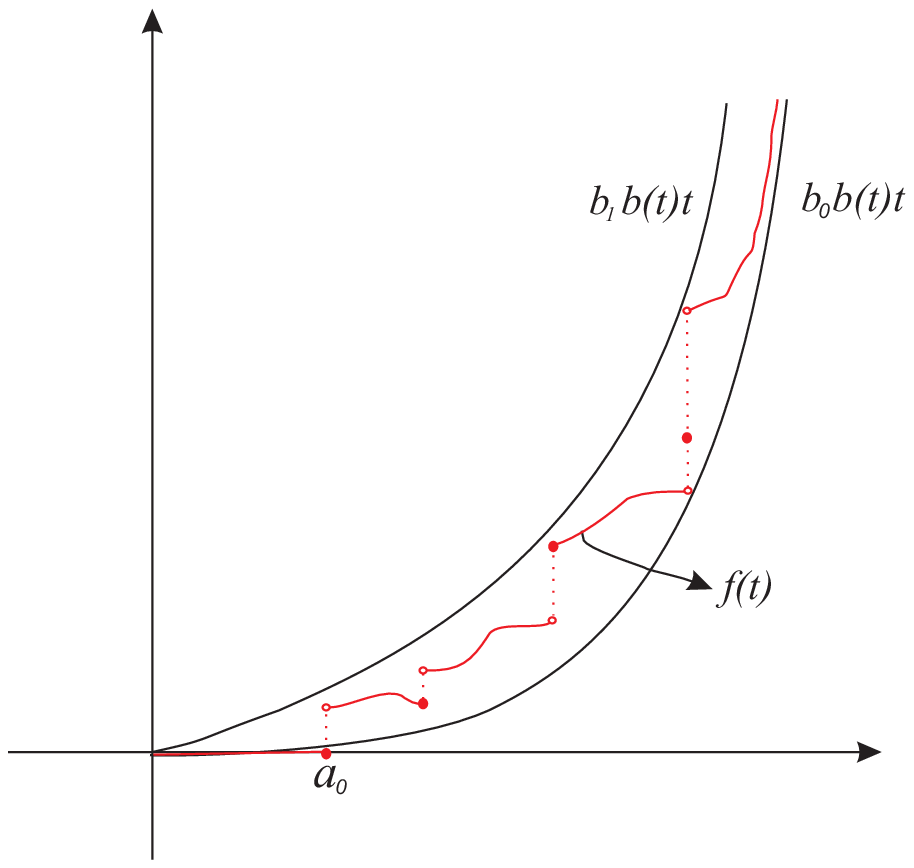}
\end{center}

In the last twenty years the study on nonlocal problems of the type
$$
\left \{ \begin{array}{l}
-M\left(\displaystyle\int_\Omega \mid\nabla u\mid^{2} dx\right)\Delta u= f(x,u)  \ \mbox{in}\ \Omega,\\
u\in H^{1}_{0}(\Omega)                                    
\end{array}
\right.\leqno{(K)}
$$
grew exponentially. That was, probably, by the difficulties existing
in this class of problems and that do not appear in the study of
local problems, as well as due to their significance in
applications. Without hope of being thorough, we mention some
articles with multiplicity results and that are related with our
main result. We will restrict our comments to the works that have
emerged in the last four years

The problem $(K)$ was studied in \cite{GJ}. The version with
p-Laplacian operator was studied in \cite{correa}. In both cases,
the authors showed a multiplicity result using genus theory. In
\cite{Xiao} the authors showed a multiplicity result
for the problem $(K)$ using the Fountain theorem and the Symmetric
Mountain Pass theorem. In all these articles the nonlinearity is
continuous. The case discontinuous was studied in \cite{Nascimento}.
With a nonlinearity of the Heaviside type the authors showed a
existence of two solutions via Mountain Pass Theorem and Ekeland's
Variational Principle.

In this work we extend the studies found in the papers above in the
following sense:

\noindent a) We cannot use the classical Clark's Theorem for $C^1$
functional (see \cite[Theorem 3.6]{Davi}), because in our case, the
energy functional is only locally Lipschitz continuous. Thus, in all
section \ref{finalissimo} we adapt for nondifferentiable functionals
an argument found in \cite{GP}.\\

\noindent b) Unlike \cite{Nascimento}, we show a result of
multiplicity using genus theory considering a nonlinearity that can
have a number enumerable of discontinuities.\\

\noindent c) Problem $(P)$ possesses more complicated
nonlinearities, for example:

\noindent (i) $\Phi(t)=t^{p_{0}}+ t^{p_{1}}$, $1< p_{0} < p_{1}< N$
and $ p_1\in (p_0,p^{*}_{0})$.

\noindent (ii) $\Phi(t)=(1+t^{2})^{\gamma}-1$, $\gamma\in
(1,\frac{N}{N-2})$.

\noindent (iii) $\Phi(t)=t^{p}\log(1+t)$ with $1<p_0<p<N-1$, where
$p_0=\frac{-1+\sqrt{1+4N}}{2}$.

\noindent (iv) $\Phi(t)=\int_0^ts^{1-\alpha}\left(\sinh^{-1}s\right)^\beta ds, \ 0\leq\alpha\leq 1, \ \beta>0.$\\

\noindent d) We work with Orlicz-Sobolev spaces and some different
estimates from those found in the papers above are necessary. For
example, the Lemma \ref{LipLoc1} is a version for Orlicz-Sobolev
spaces of a well-known result of Chang (see \cite{chang},
\cite{clarke} and \cite[Lemma 3.3]{Nascimento}). In the Lemma
\ref{PS1} one different estimate was necessary because of the
presence of the nonlocal term.

The paper is organized as follows. In the next section we present a
brief review on Orlicz-Sobolev spaces. In section  \ref{Section
results for the discontinuous} we recall some definitions  and basic
results on the critical point theory of locally Lipschitz continuous
functionals. We also present variational tools which we will prove
the main result of this paper. Furthermore, in this chapter, we
prove the Lemma \ref{LipLoc1}, which is a version for Orlicz-Sobolev
spaces of a well-known result of Chang (see \cite{chang},
\cite{clarke} and \cite[Lemma 3.3]{Nascimento}). In Section
\ref{Section Genus} we present just some preliminary results
involving genus theory that will be used in this work. In the
Section \ref{finalissimo} we prove Theorem \ref{teorema1}.


\section{A brief review on Orlicz-Sobolev spaces}\label{Section
Orlicz}

Let $\phi$ be a real-valued function defined $[0,\infty)$ and having
the following properties:

\noindent $a)$ \ \  $\phi(0)=0$, $\phi(t)>0$ if $t>0$ and
$\displaystyle\lim_{t\rightarrow \infty}\phi(t)=\infty$.

\noindent $b)$ \ \  $\phi$ is nondecreasing, that is, $s>t$ implies
$\phi(s) \geq \phi(t)$.

\noindent $c)$ \ \ $\phi$ is right continuous, that is,
$\displaystyle\lim_{s\rightarrow t^{+}}\phi(s)=\phi(t)$.

Then, the real-valued function $\Phi$ defined on $\mathbb{R}$ by
$$
\Phi(t)= \displaystyle\int^{|t|}_{0}\phi(s) \ ds
$$
is called an N-function. For an N-function $\Phi$ and  an open set
$\Omega \subseteq \mathbb{R}^{N}$, the Orlicz space
$L_{\Phi}(\Omega)$ is defined (see \cite{adams}). When $\Phi$
satisfies $\Delta_{2}$-condition, that is, when there are $t_{0}\geq
0$ and $K>0$ such that $\Phi(2t)\leq K\Phi(t)$, for all $t\geq
t_{0}$,  the space $L_{\Phi}(\Omega)$ is the vectorial space of the
measurable functions $u: \Omega \to \mathbb{R}$ such that
$$
\displaystyle\int_{\Omega}\Phi(|u|) \ dx < \infty.
$$
The space $L_{\Phi}(\Omega)$ endowed with Luxemburg norm, that is,
the norm given by
$$
|u|_{\Phi}= \inf \biggl\{\lambda >0:
\int_{\Omega}\Phi\Big(\frac{|u|}{\lambda}\Big)\ dx\leq 1\biggl\},
$$
is a Banach space. The complement function of $\Phi$, denoted by
$\widetilde{\Phi}$, is given by the Legendre transformation, that is
$$
\widetilde{\Phi}(s)=\displaystyle\max_{t \geq 0}\{st -\Phi(t)\} \ \
\mbox{for} \ \ s \geq 0.
$$
These $\Phi$ and $\widetilde{\Phi}$ are complementary each other.
Involving the functions $\Phi$ and $\widetilde{\Phi}$, we have the
Young's inequality given by
$$
st \leq \Phi(t) + \widetilde{\Phi}(s).
$$
Using the above inequality, it is possible to prove the following
H\"{o}lder type inequality
$$
\biggl|\displaystyle\int_{\Omega}u v \ dx \biggl|\leq
2|u|_{\Phi}|v|_{\widetilde{\Phi}}\,\,\, \forall \,\, u \in
L_{\Phi}(\Omega) \,\,\, \mbox{and} \,\,\,  v \in
L_{\widetilde{\Phi}}(\Omega).
$$

Hereafter, we denote by $W^{1}_{0}L_{\Phi}(\Omega)$ the
Orlicz-Sobolev space obtained by the completion of
$C^{\infty}_{0}(\Omega)$ with norm
$$
\|u\|_\Phi=|u|_\Phi+|\nabla u|_\Phi.
$$

When $\Omega$ is bounded, there is $c>0$ such that
$$
|u|_\Phi \leq c |\nabla u|_\Phi.
$$

In this case, we can consider

$$
\|u\|_{\Phi}=|\nabla u|_{\Phi}.
$$

Another important function related to function $\Phi$, is the
Sobolev conjugate function $\Phi_{*}$ of $\Phi$ defined by
$$
\Phi^{-1}_{*}(t)=\displaystyle\int^{t}_{0}\displaystyle\frac{\Phi^{-1}(s)}{s^{(N+1)/N}}ds, \ t>0.
$$
The function $\Phi_{*}$ is very important because it is related to
some embedding involving $W^{1}_{0}L_{\Phi}(\Omega)$.

We say that $\Psi$ increases essentially more slowly than $\Phi_{*}$
near infinity when
$$
\displaystyle \lim_{t \rightarrow
\infty}\frac{\Psi(kt)}{\Phi_*(t)}=0, \ \text{ for all } k>0.
$$

Let $\Omega$ be a smooth bounded domain of $\mathbb{R}^{N}$. If
$\Psi$ is any N-function increasing essentially more slowly than
$\Phi_{*}$ near infinity, then the imbedding
$W_0^{1}L_{\Phi}(\Omega)\hookrightarrow L_\Psi(\Omega)$ exists and
is compact (see \cite{adams}).

The hypotheses $(\phi_{1})-(\phi_{2})$ implies that $\Phi$,
$\widetilde{\Phi}$, $\Phi_{*}$ and $\widetilde{\Phi}_{*}$ satisfy
$\Delta_{2}$-condition. This condition allows us conclude that:

\noindent 1) $u_{n}\rightarrow 0$ in $L_{\Phi}(\Omega)$ if, and only
if, $\displaystyle\int_{\Omega}\Phi(u_{n})\ dx\rightarrow 0$.

\noindent 2)  $L_{\Phi}(\Omega)$ is separable and
$\overline{C^{\infty}_{0}(\Omega)}^{|.|_{\Phi}}=L_{\Phi}(\Omega)$.

\noindent 3)  $L_{\Phi}(\Omega)$ is reflexive and its dual is
$L_{\widetilde{\Phi}}(\Omega)$(see \cite{adams}).

\vspace{.5cm}

Under assumptions $(\phi_{1})-(\phi_{2})$, some elementary
inequalities listed in the following lemmas are valid. For the
proofs, see \cite{fukagai}.

\begin{lemma}\label{desigualdadeimportantes}
Let $\xi_{0}(t)=\min\{t^{l},t^{m}\}$,
$\xi_{1}(t)=\max\{t^{l},t^{m}\}$,
$\xi_{2}(t)=\min\{t^{l^{*}},t^{m^{*}}\}$,
$\xi_{3}(t)=\max\{t^{l^{*}},t^{m^{*}}\}$, $t\geq 0$. Then
$$
\xi_{0}(\|u\|_{\Phi})\leq \displaystyle\int_{\Omega}\Phi(|\nabla u|)
\ dx \ \leq\xi_{1}(\|u\|_{\Phi}),
$$
$$
\xi_{2}(|u|_{\Phi_{*}})\leq \displaystyle\int_{\Omega}\Phi_{*}(|u|)
\ dx \ \leq\xi_{3}(|u|_{\Phi_{*}})
$$
and
$$
\Phi_{*}(t)\geq \Phi_{*}(1)\xi_{2}(t).
$$
\end{lemma}

\begin{lemma}\label{desigualdadeimportantes1}
Let $\eta_{0}(t)=\min\{t^{q_0},t^{q_1}\},
\eta_{1}(t)=\max\{t^{q_0},t^{q_1}\}, t\geq 0$. Then
$$
\eta_{0}(|u|_{B}) \leq \displaystyle\int_{\Omega}B(|u|)dx \leq
\eta_{1}(|u|_{B})
$$
and
$$
B(1)\eta_{0}(t) \leq B(t) \leq B(1)\eta_{1}(t), \ t \in \mathbb{R}.
$$
\end{lemma}

\begin{lemma}\label{DESIGUALD}
$\widetilde{\Phi}(\frac{\Phi(s)}{s})\leq \Phi(s),\ s>0.$
\end{lemma}

The next result is a version of Brezis-Lieb's Lemma
\cite{Brezis_Lieb} for Orlicz-Sobolev spaces and the proof can be
found in \cite{Gossez}.

\begin{lemma}\label{brezislieb} Let $\Omega\subset \mathbb{R}^N$ open set and $\Phi:\mathbb{R}\rightarrow [0,\infty)$
an N-function satisfies $\Delta_2-$condition. If the
complementary function $\widetilde{\Phi}$ satisfies
$\Delta_2-$condition, $(f_n)$ is bounded in $L_\Phi(\Omega)$, such
that
$$
f_n(x)\rightarrow f(x)\  \text{a.s }x\in\Omega,
$$
then
$$
f_n\rightharpoonup f \ \text{in }L_\Phi(\Omega).
$$
\end{lemma}

\begin{coro}\label{imersao}
The imbedding $W_0^{1}L_{\Phi}(\Omega)\hookrightarrow L_{B}(\Omega)$
exists and is compact.
\end{coro}
\noindent \textbf{Proof:} It is sufficiently to show that $B$
increasing essentially more slowly than $\Phi_{*}$ near infinity.
Indeed,
$$
\frac{B(kt)}{\Phi_*(t)}\leq
\frac{B(1)B(kt)}{\xi_2(t)}=B(1)k^{q_1}t^{q_1-l^*}, \ k>0.
$$
Since $q_1<l^*$, we get
$$
\displaystyle\lim_{t\rightarrow +\infty}\frac{B(kt)}{\Phi_*(t)}= 0.
$$
\hfill\rule{2mm}{2mm}


\section{Technical results on locally Lipschitz functional and variational framework}\label{Section results for the
discontinuous}

In this section, for the reader's convenience,  we recall some
definitions and basic results  on the critical point theory of
locally Lipschitz continuous functionals  as developed by Chang
\cite{chang}, Clarke \cite{clarke, Clarke} and Grossinho \& Tersian
\cite{grossinho}. \vspace{2mm}

 Let $X$ be a real Banach space. A functional $J:X \rightarrow {\mathbb{R}}$ is locally Lipschitz continuous, $J \in
Lip_{loc}(X, {\mathbb{R}})$ for short, if given $u \in X$ there is
an open neighborhood $V := V_u \subset X$ and some constant  $K =
K_V > 0$ such that
$$
\mid J(v_2) - J(v_1) \mid \leq K \parallel v_2-v_1 \parallel,~ v_i
\in V,~ i = 1,2.
$$

\nd The directional derivative of $J$ at $u$ in the direction of $v
\in X$ is defined by

$$
J^0(u;v)=\displaystyle \displaystyle \limsup_{h \to 0,~\sigma
\downarrow 0} \frac{J(u+h+ \sigma v)-I(u+h)}{\sigma}.
$$
The generalized gradient of $J$ at $u$ is the set
$$
\partial J(u)=\big\{\mu\in X^*; \langle \mu,v\rangle\leq J^0(u;v), \ v\in
X\big \}.
$$
\nd  Since $J^0(u;0) = 0$, $\partial J(u)$ is the subdifferential of
$J^0(u;0)$. Moreover, $J^0(u;v)$ is the support function of
$\partial J(u)$ because
$$
J^0(u;v)=\max\{\langle
\xi,v\rangle; \xi\in \partial J(u)\}.
$$

The generalized gradient $\partial J(u)\subset X^*$ is convex,
non-empty and weak*-compact, and
$$
m^{J}(u) = \min \big\{\parallel\mu\parallel_{X^*};\mu \in
\partial J(u) \big \}.
$$
Moreover,
$$
\partial J(u) = \big \{J'(u) \big \},  \mbox{if}\ J \in C^1(X,{\mathbb{R}}).
$$
A critical point of $J$ is an element $u_0 \in X$ such that $0\in
\partial J(u_0)$ and a critical value of $J$ is a real number $c$
such that $J(u_0)=c$ for some critical point $u_0 \in X$.
\vspace{2mm}

About variational framework, we say that $u \in
W^{1}_{0}L_{\Phi}(\Omega)$ is a weak solution of the problem $(P)$
if it verifies
$$
M\left(\dis\int_\Omega\Phi(\mid\nabla u\mid )\ dx
\right)\dis\int_\Omega\phi(\mid\nabla u\mid)\nabla u\nabla v \
dx-\int_\Omega\phi(u)uvdx-\int_\Omega \rho v \ dx=0,
$$
for all $v \in W^{1}_{0}L_{\Phi}(\Omega)$ and for some $\rho \in
L_{\widetilde{B}}(\Omega)$ with
$$
\underline{f}(u(x))\leq \rho(x) \leq \overline{f}(u(x)) \ \
\mbox{a.e
 in} \ \ \Omega,
$$
and moreover the set $\{x\in \Omega; \mid u\mid\geq a_0\}$ has positive measure.
Thus, weak solutions of $(P)$ are critical points of the functional

$$
J(u) =\widehat{M}\left(\ds\int_{\Omega}\Phi(|\nabla u|) dx\right)
-\int_\Omega \Phi(u)dx-
          \ds\int_{\Omega}F(u)\ dx,
$$
where $\widehat{M}(t)=\displaystyle\int^{t}_{0}M(s) ds$. In order to
use variational methods, we first derive some results related to the
Palais-Smale compactness condition for the problem $(P)$.

We say that a sequence $(u_{n})\subset W^{1}_{0}L_{\Phi}(\Omega)$ is
a Palais-Smale sequence for the locally lipschitz functional $J$
associated of problem $(P)$ if
\begin{eqnarray}\label{****}
J(u_{n})\rightarrow c \ \mbox{and} \ m^{J}(u_{n})\rightarrow 0 \
\mbox{in} \ (W^{1}_{0}L_{\Phi}(\Omega))^*,
\end{eqnarray}
where $$ c = \displaystyle\inf_{\eta \in \Gamma}
\displaystyle\max_{t \in [0,1]} J(\eta(t))>0
$$
and
$$
\Gamma := \{ \eta \in C([0,1],X) : \eta(0)=0, ~I(\eta(1)) < 0\}.
$$

If (\ref{****}) implies the existence of a subsequence $(u_{n_{j}})
\subset (u_{n})$ which converges in $W^{1}_{0}L_{\Phi}(\Omega)$, we
say that these one functionals satisfies the nonsmooth $(PS)_{c}$
condition.

Note that $J \in Lip_{loc}(W^{1}_{0}L_{\Phi}(\Omega), {\mathbb{R}})$
and from convex analysis theory, for all $ w \in \partial J(u)$,
$$
\langle w,v\rangle=M\left(\ds\int_{\Omega}\Phi(|\nabla u|) \
dx\right) \ds\int_{\Omega}\phi(|\nabla u|)\nabla u \nabla v \ dx  -
\int_\Omega \phi(u)uvdx-
          \langle\rho,v\rangle,
$$
for some $\rho \in \partial \Psi(u)$, where
$\Psi(u)=\displaystyle\int_{\Omega}F(u) dx$. We have $\Psi \in
Lip_{loc}(L_{B}(\Omega), {\mathbb{R}})$, $\partial \Psi(u) \in
L_{\widetilde{B}}(\Omega)$.

The next result is a version for Orlicz-Sobolev spaces of a
well-known result of Chang (see \cite{chang}, \cite{clarke} and
\cite[Lemma 3.3]{Nascimento}).

\begin{lemma}\label{LipLoc1}
Suppose that $M_{1}$, $(f_{2})$ and $(f_{3})$ hold. For each $u \in
L_{B}(\Omega)$, if $\rho \in
\partial \Psi(u)$, then
$$
\underline{f}(u(x))\leq \rho(x) \leq \overline{f}(u(x)) \ a.e \ x
\in \Omega,
$$
and if $a_0>0$
$$
\rho(x)=0 \ a.e \ x\in \{x\in \Omega; \mid u(x)\mid < a_0\}.
$$
\end{lemma}
\noindent \textbf{Proof:} Considering $u,v \in L_{B}(\Omega)$, from
definition
\begin{eqnarray*}
\Psi^0(u;v)&=&\displaystyle\limsup_{h\rightarrow0,t\rightarrow0^+}\frac{\Psi(u+h+tv)-\Psi(u+h)}{t}\\
&=&\displaystyle\limsup_{h\rightarrow0,t\rightarrow0^+}\frac{1}{t}\int_\Omega
\left(F(u+h+tv)-F(u+h)\right)dx.
\end{eqnarray*}

We set $(h_{n}) \subset L_{B}(\Omega)$ and $(t_{n}) \subset
\mathbb{R}_{+}$ such that $h_{n}\rightarrow 0$ in $L_{B}(\Omega)$
and $t_{n}\rightarrow 0^{+}$. Thus,
\begin{equation}\label{eq 3.1}
\Psi^0(u;v)=\displaystyle\limsup_{n\rightarrow+\infty}\int_\Omega\frac{F(u+h_n+t_nv)-F(u+h_n)}{t_n}dx.
\end{equation}

Note that from the Mean Value Theorem, $(M_{1})$, $(f_2)$ and
$(f_{3})$ that,
$$
F_{n}(u,v):=\frac{F(u+h_n+t_nv)-F(u+h_n)}{t_n}\leq c
b(|\theta_{n}(x)|)|\theta_{n}(x)||v|,
$$
where
$$
\theta_{n}(x) \in \biggl[\min\bigl\{u+h_{n}+t_{n}v,
u+h_{n}\bigl\},\max\bigl\{u+h_{n}+t_{n}v, u+h_{n}\bigl\} \biggl], \
\ x \in \Omega.
$$

Using monotonicity of $b(t)t$ we get
$$
|F_{n}(u,v)|\leq  c b(|u+h_{n}+t_{n}v|)|u+h_{n}+t_{n}v||v| + c
b(|u+h_{n}|)|u+h_{n}||v|.
$$

On the other hand, by lemma \ref{DESIGUALD} we have
$$
\widetilde{B}(b(|u+h_{n}+t_{n}v|)|u+h_{n}+t_{n}v|) \leq C
B(|u+h_{n}+t_{n}v|) \leq
C\bigl(B(u)+B(h_{n})+\eta_{1}(t_{n})B(v)\bigl)
$$
and
$$
\widetilde{B}(b(|u+h_{n}+t_{n}v|)|u+h_{n}+t_{n}v|)\rightarrow
\widetilde{B}(b(|u|)|u|) \ \ a.e \ \ \mbox{in} \ \ \Omega,
$$
where $\widetilde{B}(b(|u+h_{n}+t_{n}v|)|u+h_{n}+t_{n}v| \leq c B(u)
\in L^{1}(\Omega)$.

By Lebesgue's Theorem we obtain

$$
\displaystyle\int_{\Omega}\widetilde{B}(b(|u+h_{n}+t_{n}v|)|u+h_{n}+t_{n}v|)
dx \rightarrow \displaystyle\int_{\Omega}\widetilde{B}(b(|u|)|u|)
dx.
$$

From \ref{brezislieb} we conclude that

$$
\displaystyle\int_{\Omega}\widetilde{B}\left(b(|u+h_{n}+t_{n}v|)|u+h_{n}+t_{n}v|
- b(|u|)|u| \right) dx \rightarrow 0.
$$

Moreover, with obvious changes, we can prove that

$$
\displaystyle\int_{\Omega}\widetilde{B}\left(b(|u+h_{n}|)|u+h_{n}| -
b(|u|)|u| \right) dx \rightarrow 0.
$$

Thus, by Fatou's lemma that
 \begin{eqnarray}\label{G}
 \displaystyle\limsup
\int_{\Omega} F_{n}(u,v) \ dx \leq \displaystyle\int_{\Omega}
\limsup F_{n}(u,v) \ dx.
\end{eqnarray}

From (\ref{eq 3.1}) and (\ref{G}) we get
$$
\Psi^{0}(u,v) \leq \displaystyle\int_{\Omega} F^{0}(u,v) \ dx =
\displaystyle\int_{\Omega}\max \{\langle\xi,v\rangle; \xi \in
\partial F(u)\}dx.
$$

Consider $\widehat{\rho} \in \partial \Psi(u) \subset L_{B}(\Omega)^* \equiv
L_{\widetilde{B}}(\Omega)$ with $u \in L_{B}(\Omega)$. Then, there
is $\rho \in L_{\widetilde{B}}(\Omega)$ such that
$$
\langle\widehat{\rho}, v\rangle= \displaystyle\int_{\Omega} \rho v \
dx, \ \ v
 \in  L_{B}(\Omega).
$$

We claim that
$$
\rho(x) \geq \underline{f}(u(x)) \ \ a.e \ \ \mbox{in} \ \ \Omega.
$$

Arguing, by contradiction, we suppose that there is $A \subset
\Omega$ with $|A|>0$ such that $\rho(x) < \underline{f}(u(x))$.
Hence,
\begin{eqnarray}\label{G1}
\displaystyle\int_{A}\rho(x) \ dx<
\displaystyle\int_{A}\underline{f}(u(x)) \ dx.
\end{eqnarray}

Let $v=-\chi_{A}$ be a function in $ L_{B}(\Omega)$, where
$\chi_{A}$ is characteristic  function of set $A$. Thus,
$$
-\displaystyle\int_{A}\rho \ dx =\displaystyle\int_{\Omega}\rho v \
dx \leq \Psi^{0}(u,v) \leq
\displaystyle\int_{\Omega}\underline{f}(u(x)) v dx = -
\displaystyle\int_{A}\underline{f}(u(x)) \ dx,
$$
with is a contradiction with $(\ref{G1})$. Thus
$$
\rho(x) \geq \underline{f}(u(x)) \ \ a.e \ \ \mbox{in} \ \ \Omega.
$$

The inequality
$$
\rho(x) \leq \overline{f}(u(x)) \ \ a.e \ \ \mbox{in} \ \ \Omega
$$
follows the same argument.{\fim}


\section{Results involving genus}\label{Section Genus}

We will start by considering some basic notions on the Krasnoselskii
genus that we will use in the proof of our main results.

Let $E$ be a real Banach space. Let us denote by $\mathfrak{A}$ the
class of all closed subsets  $A\subset E\setminus \{0\}$ that are
symmetric with respect to the origin, that is, $u\in A$ implies
$-u\in A$.

\begin{definition}
Let $A\in \mathfrak{A}$. The Krasnoselskii genus $\gamma(A)$ of $A$
is defined as being the least positive integer $k$ such that there
is an odd mapping $\phi \in C(A,\mathbb{R}^{k})$ such that
$\phi(x)\neq 0$ for all $x\in A$. If $k$ does not exist we set
$\gamma(A)=\infty$. Furthermore, by definition,
$\gamma(\emptyset)=0$.
\end{definition}

In the sequel we will establish only the properties of the genus
that will be used through this work. More information on this
subject may be found in the references  by  \cite{Ambrosetti},
\cite{Castro}, \cite{Davi} and \cite{Kranolseskii}.
\begin{prop}
Let $E={\mathbb{R}}^{N}$ and $\partial\Omega$ be the boundary of an
open, symmetric and bounded subset $\Omega \subset {\mathbb{R}}^{N}$
with $0 \in \Omega$. Then $\gamma(\partial\Omega)=N$.
\end{prop}

\begin{coro}\label{esfera}
$\gamma(\mathcal{S}^{N-1})=N$ where $\mathcal{S}^{N-1}$ is a unit
sphere of ${\mathbb{R}}^{N}$.
\end{coro}

\begin{prop}\label{paracompletar}
If $K \in \mathfrak{A}$, $0 \notin K$ and $\gamma(K) \geq 2$, then
$K$ has infinitely many points.
\end{prop}


\section{Proof of Theorem \ref{teorema1}}\label{finalissimo}


The plan of the proof is to show that the set of critical points of
the functional $J$ is compact, symmetric, does not contain the zero
and has genus  more than $2$. Thus, our main result is a consequence
of Proposition \ref{paracompletar}.

In the proof of the Theorem  \ref{teorema1} we shall need the
followings technical results:

\begin{lemma}\label{Gl}
The functional $J$ is coercive.
\end{lemma}
\noindent \textbf{Proof:} Using $(M_{1})$ and $(f_{2})$  we get
\begin{eqnarray*}
 J(u)&\geq& k_0\int_0^{\displaystyle\int_\Omega\Phi(\mid\nabla
u\mid) \ dx}s^\alpha ds-\int_\Omega\Phi(u)dx- b_1\displaystyle\int_\Omega B(u)\ dx\\
&\geq&\frac{k_0}{\alpha+1}\left(\int_\Omega \Phi(\mid\nabla
u\mid)dx\right)^{\alpha +1}-\int_\Omega\Phi(u)dx-b_1\int_\Omega B(u)
\ dx .
 \end{eqnarray*}

From Lemmas \ref{desigualdadeimportantes} and
\ref{desigualdadeimportantes1} we obtain
$$
J(u)\geq \frac{k_0}{\alpha
+1}\xi_0(\|u\|_\Phi)^{\alpha+1}-\xi_1(\mid u\mid_\Phi)-b_1\eta_1(|
u|_B).
$$

Using now Corollary \ref{imersao} we get the continuous imbedding
$W_0^{1}L_{\Phi}(\Omega)\hookrightarrow
L_{B}(\Omega),L_\Phi(\Omega)$ hold. Hence, there are positive
constants $C_{1},C_{2}$ and $C_{3}$ such that, for $|\nabla
u|_{\Phi}\geq 1$, we have
\begin{eqnarray*}
J(u)&\geq& C_{1}\| u\|_{\Phi}^{l(\alpha+1)}-C_{2}\| u\|_{\Phi}^m-
C_3\| u\|_{\Phi}^{q_1}.
\end{eqnarray*}
Since  $l(\alpha+1)>q_1>m$, we conclude that $J$ is
coercive.\hfill\rule{2mm}{2mm}

Now we prove that $J$ satisfies the nonsmooth $(PS)_{c}$ condition.

\begin{lemma}\label{PS1}
The functional $J$ satisfies the nonsmooth $(PS)_{c}$ condition,
for all $c \in \mathbb{R}$.
\end{lemma}
\noindent \textbf{Proof:} Let $(u_{n})$ be a sequence in
$W^{1}_{0}L_{\Phi}(\Omega)$ such that
$$
J(u_{n})\rightarrow c \ \ \mbox{and} \ \ m^{J}(u_{n})\rightarrow 0.
$$

From now we consider $(w_{n})\subset\partial J(u_n)\subset
(W^{1}_{0}L_{\Phi}(\Omega))^{*}$ such that
$$
m^{J}(u_{n})=\|w_{n}\|_{*}=o_{n}(1)
$$
and
$$
\langle w_{n},v\rangle=
M\left(\displaystyle\int_{\Omega}\Phi(|\nabla u_{n}|) \
dx\right)\displaystyle\int_{\Omega}\phi(|\nabla u_{n}|)\nabla
u_{n}\nabla v \ dx-\int_\Omega \phi(u)uvdx
-\langle\rho_{n},v\rangle,
$$
with $\rho_{n} \in \partial\Psi(u_{n})$.

Note that from Lemma \ref{LipLoc1} we have
$$
\underline{f}(u_{n}) \leq \rho_{n} \leq \overline{f}(u_{n}) \ \
\mbox{a.e in} \ \ \Omega.
$$

On the other hand, since $J$ is coercive, we derive that $(u_{n})$
is bounded in $W^{1}_{0}L_{\Phi}(\Omega)$. Thus, passing to a
subsequence, if necessary, we have
$$
u_{n}\rightharpoonup u \ \ \mbox{in} \ \ W^{1}_{0}L_{\Phi}(\Omega),
$$
$$
\frac{\partial u_{n}}{\partial x_{i}}\rightharpoonup \frac{\partial
u}{\partial x_{i}} \ \ \mbox{in} \ \ L_{\Phi}(\Omega),
$$
$$
u_{n}\rightarrow u \ \ \mbox{in} \ \ L_{B}(\Omega) \mbox{ and } L_\Phi(\Omega),
$$
and
$$
\displaystyle\int_{\Omega}\Phi (|\nabla u_{n}|) \ dx \rightarrow
t_{0} \geq 0.
$$

If $t_{0}=0$, then from Lemma \ref{desigualdadeimportantes} we
obtain
$$
| \nabla u_n|_\Phi\leq \xi_0^{-1}\left(\displaystyle\int_\Omega
\Phi(\mid\nabla u_n\mid)  \ dx \right) \rightarrow 0
$$
and the proof is finished.

If $t_{0}>0$, since $M$ is a continuous function, we get
$$
M\left(\displaystyle\int_\Omega\Phi(\mid\nabla u_n\mid) \ dx
\right)\rightarrow M(t_0).
$$

Thus, from $(M_{1})$ and for $n$ sufficiently large,
\begin{eqnarray}\label{limitacaoporbaixo1}
M\left(\int_\Omega\Phi (\mid\nabla u_n\mid) \ dx \right)\geq
k_0t_{0}^{\alpha}>0.
\end{eqnarray}

Now we proof that $(\rho_{n})$ is bounded in
$L_{\widetilde{B}}(\Omega)$. Note that from $(f_{2})$, $(f_{3})$ and
$(M_{1})$ that
$$
\overline{f}(t) \leq c b(t)t.
$$

Since
$$
\overline{f}(t)=-\underline{f}(-t)
$$
we get from Lemmas \ref{desigualdadeimportantes1} and \ref{DESIGUALD} that
\begin{eqnarray*}
\displaystyle\int_{\Omega}\widetilde{B}(u_{n}) \ dx &\leq & C
\displaystyle\int_{[u_{n}\geq 0]}B(u_{n}) \ dx+
\displaystyle\int_{[u_{n}< 0]}\widetilde{B}(\underline{f}(- u_{n}))
\ dx \\
&\leq & \displaystyle\int_{\Omega}\widetilde{B}(u_{n}) \ dx \leq
C\eta (|u_{n}|_{B}) \leq \overline{C}(\|u_{n}\|_{\Phi}),
\end{eqnarray*}
which implies that $(\rho_{n})$ is bounded in
$L_{\widetilde{B}}(\Omega)$.

Then
\begin{eqnarray}\label{convergrho}
\displaystyle\int_{\Omega}\rho_{n}(u_{n}-u) \ dx \rightarrow 0.
\end{eqnarray}

From definition of $(u_{n})$ we have
\begin{eqnarray*}
o_{n}(1)= \langle w_{n}, u_{n}-u\rangle&= &
M\left(\displaystyle\int_{\Omega}\Phi(|\nabla u_{n}|) \
dx\right)\displaystyle\int_{\Omega}\phi(|\nabla u_{n}|)\nabla
u_{n}\nabla (u_{n}- u) \ dx\nonumber\\
&&-\int_\Omega\phi(u_n)u_n(u_n-u)dx- \displaystyle\int_{\Omega}\rho_{n} (u_{n}- u) \ dx.
\end{eqnarray*}
Since $|u_n-u|_\Phi$ goes to $0$ and $(\phi(u_n)u_n)$ is bounded in
$L_{\widetilde{\Phi}}(\Omega)$, have that
\begin{eqnarray}\label{der01}
\int_\Omega\phi(u_n)u_n(u_n-u)dx\rightarrow 0.
\end{eqnarray}
We get from  (\ref{convergrho}) and (\ref{der01}) that
$$
M\left(\displaystyle\int_{\Omega}\Phi(|\nabla u_{n}|) \
dx\right)\displaystyle\int_{\Omega}\phi(|\nabla u_{n}|)\nabla
u_{n}\nabla (u_{n}- u) \ dx\rightarrow 0.
$$

From (\ref{limitacaoporbaixo1}) and the last convergence implies
that
$$
\displaystyle\int_{\Omega}\phi(|\nabla u_{n}|)\nabla u_{n}\nabla
(u_{n}- u) \ dx\rightarrow 0.
$$
Setting $\beta:{\mathbb{R}}^N\rightarrow{\mathbb{R}} ^N$ by
$$
\beta(x)=\phi(\mid \nabla x\mid)\nabla x, \ x\in{\mathbb{R}}^N,
$$
the last limit imply that for some subsequence, still denoted by
itself,
$$
\left(\beta(\nabla u_n(x))-\beta(\nabla u(x))\right)(\nabla
u_n(x)-\nabla u(x))\to 0 \ \mbox{a.e in  } \,\,\, \Omega.
$$
Applying a result found in Dal Maso and Murat \cite{Maso}, it
follows that
$$
\nabla u_n(x)\to \nabla u(x) \ \mbox{a.e in } \,\,\,\Omega.
$$
Then
$$
u_{n}\rightarrow u \ \ \mbox{in} \ \ W^{1}_{0}L_{\Phi}(\Omega).
$$
\hfill\rule{2mm}{2mm}

Let $K_{c}$ be the set of critical points of $J$. More precisely
$$
K_{c}=\{u \in W^{1}_{0}L_{\Phi}(\Omega): 0 \in \partial J(u) \ \
\mbox{and} \ \ J(u)=c\}.
$$

Since $J$ is even, we have that $K_{c}$ is symmetric. The next
result is important in our arguments and allows  we conclude that
$K_{c}$ is compact. The proof can be found in \cite{chang}.
\begin{lemma}\label{KcCompacto}
If $J$ satisfies the nonsmooth $(PS)_{c}$ condition, then $K_{c}$ is
compact.
\end{lemma}

To prove that $K_{c}$ does not contain zero, we construct a special
class of the levels $c$.

For each $k \in \mathbb{N}$, we define the set
$$
\Gamma_{k}=\{C \subset W^{1}_{0}L_{\Phi}(\Omega): C \ \ \mbox{is
closed}, C=-C \ \ \mbox{and} \ \ \gamma(C) \geq k\},
$$
and the values
$$
c_{k}=\displaystyle\inf_{C\in \Gamma_{k}}\displaystyle\sup_{u \in
C}J(u).
$$

Note that
$$
-\infty\leq c_{1}\leq c_{2}\leq c_{3}\leq ...\leq
c_{k}\leq ...
$$
and, once that $J$ is coercive and continuous, $J$ is bounded below
and, hence, $c_{1} > -\infty$. In this case, arguing as in
\cite[Proposition 3.1]{BWW}, we can prove that each $c_{k}$ is a
critical value for the functional $J$.

\begin{lemma}\label{minimax}
Given $k \in \mathbb{N}$, there exists $\epsilon = \epsilon(k)>0$
such that
$$
\gamma(J^{-\epsilon}) \geq k,
$$
where $J^{-\epsilon}=\{u \in W^{1}_{0}L_{\Phi}(\Omega): J(u) \leq
-\epsilon\}$.
\end{lemma}
\noindent\textbf{Proof:} Fix $k \in \mathbb{N}$, let $X_{k}$ be a
k-dimensional subspace of $W^{1}_{0}L_{\Phi}(\Omega)$. Thus, there
exists $C_k>0$ such that
$$
-C_k|\nabla u|_{\Phi}\geq - \displaystyle|u|_{\Phi},
$$
for all $u \in X_{k}$.

We now use the inequality above, $(M_{1})$, $(f_{2})$, $(f_{3})$,
Lemmas \ref{desigualdadeimportantes} and
\ref{desigualdadeimportantes1}  to conclude that
\begin{eqnarray*}
J(u)\leq \frac{k_1}{\alpha+1}\xi_1(\|
u\|_\Phi)^{\alpha+1}-\xi_0(C_k\|u\|_\Phi).
\end{eqnarray*}

For $\|u\|_{\Phi} \leq 1$ we get
$$
J(u)\leq \| u \|_\Phi^{m}\left(\frac{k_1}{\alpha+1}\|
u\|_\Phi^{(\alpha+1)l-m}-C_k^{m}\right).
$$
Considering $R>0$ such that
$$
R<min\left\{1,
\left(\frac{\alpha+1}{k_1}C_k^{m}\right)^{\frac{1}{(\alpha+1)l-m}}\right\},
$$
there exists $\epsilon=\epsilon(R)>0$ such that
$$
J(u)<-\epsilon < 0,
$$
for all $u\in {\mathcal{S}_R}=\{u\in X_k; |\nabla u |_\Phi=R \}$.
Since $X_k$ and $\mathbb{R}^k$ are isomorphic and $\mathcal{S}_R$
and $S^{k-1}$ are homeomorphic, we conclude from Corollary
\ref{esfera} that $\gamma(\mathcal{S}_R)=\gamma(S^{k-1})=k$.
Moreover, once that ${\mathcal{S}_R} \subset J^{-\epsilon}$ and
$J^{-\epsilon}$ is symmetric and closed,  we have
$$
k= \gamma ({\mathcal{S}_R})\leq \gamma( J^{-\epsilon}).
$$
\hfill\rule{2mm}{2mm}

\begin{lemma}\label{minimax1}
Given $k \in \mathbb{N}$, the number $c_{k}$ is negative.
\end{lemma}
\noindent\textbf{Proof:} From Lemma \ref{minimax}, for each $k\in
\mathbb{N}$ there exists $\epsilon >0$ such that
$\gamma(J^{-\epsilon}) \geq k$. Moreover, $ 0 \notin J^{-\epsilon}$
and $J^{-\epsilon}\in \Gamma_{k}$. On the other hand
$$
\displaystyle\sup_{u\in J^{-\epsilon}}J(u)\leq -\epsilon.
$$

Hence,
$$
-\infty < c_{k}=\displaystyle\inf_{C\in
\Gamma_{k}}\displaystyle\sup_{u \in C}J(u) \leq
\displaystyle\sup_{u\in J^{-\epsilon}}J(u) \leq -\epsilon <0.
$$
\hfill\rule{2mm}{2mm}

A direct consequence of the last Lemma is that $0 \notin K_{c_{k}}$.
The next result is also important in our arguments and the proof can
be found in \cite{chang}.

\begin{lemma}\label{Deformacao}
Suppose that $X$ is a reflexive Banach space and $J$ is even and a
locally Lipschitz function, satisfying the $(PS)_{c}$ condition. If
$U$ is any neighborhood of $K_{c}$, then for any $\epsilon_0>0$
there exist $\epsilon \in (0,\epsilon_{0})$ and a odd homeomorphism
$\eta :
X\rightarrow X$ such that:\\
\noindent $a)$ $\eta (x)=x$ for $x\notin
J^{c+\epsilon}\backslash J^{c-\epsilon}$ \\
\noindent $b)$ $\eta(J^{c+\epsilon} \backslash U)\subset J^{c-\epsilon}$\\
\noindent $c)$ If $K_{c}=\emptyset$, then
$\eta(J^{c+\epsilon})\subset J^{c-\epsilon}$.
\end{lemma}

\begin{lemma}\label{minimax2}
If $c_{k}=c_{k+1}=...=c_{k+r}$ for some $r \in \mathbb{N}$, then
$$
\gamma(K_{c_{k}})\geq r+1.
$$
\end{lemma}
\noindent\textbf{Proof:} Suppose, by contradiction, that
$\gamma(K_{c_k})\leq r$. Since $K_{c_{k}}$ is compact and symmetric,
there exists a closed and symmetric set $U$ with $ K_{c_{k}}\subset
U$ such that $\gamma(U)= \gamma(K_{c_{k}}) \leq r$. Note that we can
choose $U\subset J^{0}$ because $c_{k}<0$. By the deformation lemma
\ref{Deformacao} we have an odd homeomorphism $ \eta:
W^{1}_{0}L_{\Phi}(\Omega)\rightarrow W^{1}_{0}L_{\Phi}(\Omega)$ such
that $\eta(J^{c_{k}+\delta}-U)\subset J^{c_{k}-\delta}$ for some
$\delta
> 0$ with $0<\delta < -c_{k}$. Thus, $J^{c_{k}+\delta}\subset J^{0}$ and by
definition of $c_{k}=c_{k+r}$, there exists $A \in \Gamma_{k+r}$
such that $\displaystyle\sup_{u \in A}J(u) < c_{k}+\delta$, that is,
$A \subset J^{c_{k}+\delta}$ and
\begin{eqnarray}\label{estrela1}
\eta(A-U) \subset \eta ( J^{c_{k}+\delta}-U)\subset
J^{c_{k}-\delta}.
\end{eqnarray}
But $\gamma(\overline{A-U})\geq \gamma(A)-\gamma(U) \geq k$ and
$\gamma(\eta(\overline{A-U}))\geq  \gamma(\overline{A-U})\geq k$.
Then $\eta(\overline{A-U}) \in \Gamma_{k}$ and this contradicts
(\ref{estrela1}). Hence, this lemma is proved.\hfill\rule{2mm}{2mm}

\subsection{Proof of Theorem \ref{teorema1}}

If $-\infty< c_{1} < c_{2} < ...< c_{k}< ...<0$ and since each
$c_{k}$ critical value of $J$, then we obtain infinitely many
critical points of $J$ and hence, the problem $(P)$ has infinitely
many solutions.

On the other hand, if there are two constants $c_{k}=c_{k+r}$, then
$c_{k}=c_{k+1}=...=c_{k+r}$ and from Lemma \ref{minimax2}, we have
$$
\gamma(K_{c_{k}})\geq r+1 \geq 2.
$$
From Proposition \ref{paracompletar}, $K_{c_{k}}$ has infinitely
many points.

Let $(u_{k})$ critical points of $J$. Now we show that, for
\begin{eqnarray}\label{setmeasure}
a_0 < \xi_1^{-1}\left(\frac{k_0 l}{\Phi(1)\mid \Omega\mid
m^2}\xi_0(C)^{\alpha+1}\right),
\end{eqnarray}
we have that
$$
\bigl\{x \in \Omega: \mid u_{k}(x)\mid \geq a_0\bigl\}
$$
has positive measure. Thus every critical points of $J$, are solutions of $(P)$. Suppose, by contradiction, that this set has
null measure. Thus
\begin{eqnarray*}
0&=&M\left(\int_\Omega\Phi(\mid \nabla u_k\mid )dx\right)\int_\Omega
\phi(\mid\nabla u_k\mid)\mid\nabla u_k\mid^2dx-\int_\Omega \phi(\mid
u_k\mid)\mid
u_k\mid^2 dx \\
  &\geq &  k_0 l\left(\int_\Omega\Phi(\mid\nabla
u_k\mid)dx\right)^{\alpha+1}-m\int_\Omega \Phi(u_k)dx\\
 &\geq & k_0 l\xi_0(\parallel
u_k\parallel_\Phi)^{(\alpha+1)}-m\Phi(a_0)\mid \Omega\mid,\
 \end{eqnarray*}
where we conclude
\begin{eqnarray}\label{setmeasure1}
k_0 l\xi_0(\parallel u_k\parallel_\Phi)^{(\alpha+1)}\leq
m\xi_1(a_0)\mid \Omega\mid\Phi(1).
\end{eqnarray}

Since $c_{k}\leq -\epsilon<0$, there exists $C>0$ such that
$\|u_{k}\|\geq C>0$. Hence
\begin{eqnarray*}
a_0 \geq \xi_1^{-1}\left(\frac{k_0 l}{\Phi(1)\mid \Omega\mid
m^2}\xi_0(C)^{\alpha+1}\right),
\end{eqnarray*}
which contradicts (\ref{setmeasure}). Then, $$ \bigl\{x \in \Omega:
\mid u_{k}(x)\mid \geq a_0\bigl\}
$$
has positive measure. \hfill\rule{2mm}{2mm}

\end{document}